\documentclass[journal]{IEEEtran}

\usepackage{cite}
\usepackage{graphicx}
\usepackage{color}
\DeclareGraphicsExtensions{.pdf}
\usepackage{fixltx2e}
\usepackage{url}
\usepackage{eurosym}
\usepackage[utf8]{inputenc}
\DeclareUnicodeCharacter{20AC}{\euro}
\usepackage{booktabs}
\usepackage{multirow}
\usepackage[table]{xcolor}
\usepackage{latexsym}
\usepackage{amsmath}
\usepackage{siunitx} 
\usepackage{mathtools}

\usepackage{enumitem}
\setlist[description]{leftmargin=0cm,labelindent=0.6cm}

\usepackage{multirow}
\usepackage{amssymb}
\usepackage{array}
\newcolumntype{L}[1]{>{\raggedright\let\newline\\\arraybackslash\hspace{0pt}}m{#1}}
\newcolumntype{C}[1]{>{\centering\let\newline\\\arraybackslash\hspace{0pt}}m{#1}}
\newcolumntype{R}[1]{>{\raggedleft\let\newline\\\arraybackslash\hspace{0pt}}m{#1}}

\renewcommand{\b}{\boldsymbol }

\usepackage{algorithm}
\usepackage{algorithmic}

\makeatletter
\newcommand\fs@norules{\def\@fs@cfont{\bfseries}\let\@fs@capt\floatc@ruled
  \def\@fs@pre{}%
  \def\@fs@post{}%
  \def\@fs@mid{\kern3pt}%
  \let\@fs@iftopcapt\iftrue}
\makeatother
\floatstyle{norules}
\restylefloat{algorithm}

\begin{document}
%
\title{Short-term Forecasting of Price-responsive Loads Using Inverse Optimization}
%
%
%

\author{Javier~Saez-Gallego,
        Juan~M.~Morales
\thanks{J. Saez-Gallego, J. M. Morales (corresponding author) are with the Technical University of Denmark, DK-2800 Kgs. Lyngby, Denmark (email addresses: \{jsga, jmmgo\}@dtu.dk), and their work is partly funded by DSF (Det Strategiske Forskningsr{\aa}d) through the CITIES research center (no. 1035-00027B) and the iPower platform project (no. 10-095378).}
}

\markboth{}%
{Saez-Gallego \MakeLowercase{\textit{et al.}}: Short-term Forecasting of Price-responsive Loads Using Inverse Optimization}

\maketitle
\begin{abstract}
We consider the problem of forecasting the aggregate demand of a pool of price-responsive consumers of electricity. The price-response of the aggregation is modeled by an optimization problem that is characterized by a set of marginal utility curves and minimum and maximum power consumption limits.
The task of estimating these parameters is addressed using a generalized inverse optimization scheme that, in turn, requires solving a nonconvex mathematical program. We introduce a solution method that overcomes the nonconvexities by solving instead two linear problems with a penalty term, which is statistically adjusted by using a cross-validation algorithm.
The proposed methodology is data-driven and leverages information from regressors, such as time and weather variables, to account for changes in the parameter estimates.
The power load of a group of heating, ventilation, and air conditioning systems in buildings is simulated, and the results show that the aggregate demand of the group can  be successfully captured by the proposed model, making it suitable for short-term forecasting purposes.
\end{abstract}



\IEEEpeerreviewmaketitle

\section*{Notation}

The notation used throughout the paper is stated below for quick reference. Other symbols are defined as required.
\subsection{Indexes}
\begin{IEEEdescription}
\item[$t$] Time period, ranging from 1 to $T$.
\item[$b$] Marginal utility block, ranging from 1 to $B$.
\item[$r$] Regressor, ranging from 1 to $R$.
\end{IEEEdescription}

\subsection{Decision variables}
\begin{IEEEdescription}
\item[$x_{b,t}$] Load from energy block $b$ and time $t$.
\item[$ \underline{P}_t $]Lower  bound for electricity consumption at time $t$.
\item[$\overline{P}_t$] Upper bound for electricity consumption at time $t$.
\item[$u_{b,t}$] Marginal utility of load block $b$ at time $t$.
\item[$\underline{\mu}$] Intercept for the lower load-consumption bound.
\item[$\overline{\mu}$] Intercept for the upper load-consumption bound.
\item[$\mu^u$] Intercept for marginal utility.
\item[$\underline{\alpha}_r$] Coefficient relative to the affine dependence of the lower load-consumption bound on regressor $r$.
\item[$\overline{\alpha}_r$] Coefficient relative to the affine dependence of the upper load-consumption bound on regressor $r$.
\item[$\alpha^u_r$] Coefficient relative to the affine dependence of marginal utility on regressor $r$.
\item[$\epsilon_t$] Duality gap at time $t$.
\item[$\underline{\lambda}_{t}$] Dual variable associated with the lower bound for total load at time $t$.
\item[$\overline{\lambda}_{t}$] Dual variable associated with the upper bound for total load at time $t$.
\item[$\underline{\phi}_{b,t}$] Dual variable associated with the positive block-size constraint for block $b$ at time $t$.
\item[$\overline{\phi}_{b,t}$] Dual variable associated with the maximum block-size constraint for block $b$ at time $t$.
\item[$\overline{\xi}^+_t$] Feasibility slack variable linked to the upper load-consumption bound at time $t$.
\item[$\underline{\xi}^+_t$] Feasibility slack variable linked to the lower load-consumption bound at time $t$.
\item[$\overline{\xi}^-_t$] Infeasibility slack variable linked to the upper load-consumption bound at time $t$.
\item[$\underline{\xi}^-_t$] Infeasibility slack variable linked to the lower load-consumption bound at time $t$.
\end{IEEEdescription}

\subsection{Parameters}
\begin{IEEEdescription}
\item[$x'_{t}$] Measured load at time $t$.
\item[$\tilde{x}_{b,t}'$] Adjusted measured load for block $b$ at time $t$.
\item[$p_{t}$] Price of electricity at time $t$.
\item[$Z_{r,t}$] Value of regressor $r$ at time $t$.
\item[$E_{b,t}$] Width of load block $b$ at time $t$.
\item[K] Feasibility penalty parameter.
\end{IEEEdescription}

\section{Introduction}



Demand response programs aim to alter the power consumption profile of end-users by external stimulus \cite{o2014benefits}, with the final goal of avoiding over-investing in transmission lines and generating capacities. A popular scheme amongst the numerous programs for demand-side management is \emph{Real-time Pricing} (RTP), where the external stimulus consists of varying prices along the day reflecting the change of balance between supply and demand \cite{vardakas2015survey,borenstein2005long}. Consumers of electricity, equipped with a smart grid meter and an Energy Management Controller (EMC), seek the most favorable pattern of consumption according to the dynamic price. In the case of households, the EMC comprises a home automation equipment that considers both the price of electricity and the personal preferences of the users to optimally schedule their electricity demand needs and their appliances \cite{mohsenian2010autonomous,ConejoMorales2010}. All in all, under the RTP paradigm, the consumers are \emph{price-responsive}.

Forecasting the expected electricity demand at aggregate levels, i.e., \emph{load forecasting}, is of utmost importance for network operators to enhance planning, for example, by mitigating grid congestion during peak-demand periods. Also, it is widely used by electric utilities to minimize the costs of over- or under-contracting power in electricity markets.
The increasing penetration of smart grid technologies call for solutions able to forecast the aggregate price-responsive load as accurately as possible.


In response to these challenges, the contributions of this paper are threefold:
\begin{enumerate}
\item A methodology to forecast the aggregate consumption of a cluster of price-responsive power loads using inverse optimization.
\item A computationally efficient method that approximates the solution to a generalized inverse optimization problem by solving instead two linear programming problems. The proposed approach relies on cross-validation techniques to optimally tune a penalty parameter so that the out-of-sample forecasting error is minimized.
\item A comprehensive analysis of the performance of the proposed forecasting methodology. The analysis is based on a case study that considers the simulated price-response of a group of buildings equipped with heat pumps. Furthermore, we benchmark our methodology against persistence forecasting and a state-of-the-art autoregressive moving average model with exogenous inputs \cite{corradi2012}.
\end{enumerate}


The presented methodology relates to the existing literature in several aspects. First of all, its final goal is to predict a demand for electricity, hence, it fits into the realm of load forecasting. Amongst the vast load forecasting literature \cite{weron2007modeling}, there are some authors that in the last years have focused on modeling the effect of the price on the load, for example, using a B-spline approach \cite{Hosking}, an Auto-Regressive Model With eXogeneous Inputs (ARX) \cite{corradi2012}, neural networks \cite{yun2008rbf}, or a hybrid approach with data association mining algorithms \cite{motamedi2012electricity}. The novelty of the proposed methodology in this paper with respect to the existing literature, lies in the characterization of the response of the load to price by an optimization problem. Indeed, to the best of our knowledge, we are the first ones to exploit inverse optimization for time series forecasting and, in particular, for load prediction.

The first formal description of the inverse linear programming problem is given by \cite{AhujaInverse}, which seeks to find the minimal perturbation of the objective function cost vector that makes a given data point optimal. More recent works address the case where the observations are noisy and an exact solution of the inverse problem might not exist \cite{Saez2015,chan2014generalized,Keshavarz_imputinga,aswani2015inverse,esfahani2015data,bertsimas2015data}. The proposed methodology neither makes any assumption on whether the data measurements are noisy or not, nor on the existence of a solution to the exact inverse optimization problem.

Here, as in \cite{Saez2015,xudata2015}, we extend the concept of inverse optimization to the case where right-hand side parameters of the forward linear programming problem are also to be estimated. Authors of \cite{xudata2015} assume that a feasible region for the forward problem exists, whereas we do not make any assumption in this regard. Indeed, we calculate the best feasible region, in terms of forecasting capabilities, even though it makes the observed data infeasible.
The novelty of this paper with respect to \cite{Saez2015} is twofold. First, we propose an inverse optimization scheme that is especially tailored to one-step ahead forecasting, and not to market bidding. Second, the estimation problem we formulate is not based on relaxing the KKT conditions of the forward problem. Instead, we statistically determine the feasible region and the objective function of the forward problem that render the best out-of-sample prediction performance.



%











The rest of this paper is structured as follows. In Section \ref{sec:method} we provide a general overview of the proposed forecasting methodology and the associated estimation problem. Then, in Section \ref{sec:applied}, the specific load forecasting model is provided. Section \ref{sec:simulation} introduces the framework we have used to simulate the price-response of a group of buildings equipped with heat pumps. In Section \ref{sec:study_case}, we discuss results from a case study, and finally, in Section \ref{sec:conclusions}, conclusions are duly drawn.

\section{Inverse Optimization Methodology} \label{sec:method}

Next we introduce the problem of forecasting using inverse optimization and describe the methodology that is applied later, in Section \ref{sec:applied}, to predict price-responsive electricity load.

We start from the premise that the choices made by a certain decision-maker (e.g., an aggregation of price-responsive power loads) at a certain time $t$, denoted by $\boldsymbol x_t$, are driven by the solution to the following linear optimization problem:
\begin{IEEEeqnarray}{lrl}
\text{RP$_t$($\rho_t|\boldsymbol c, \boldsymbol b$):} \qquad & \underset{\boldsymbol x_t}{\text{Maximize}} \ & (\boldsymbol c-\rho_t \boldsymbol e)^T \boldsymbol x_t \notag
\\
& \text{subject to} \  & \boldsymbol A \boldsymbol x_t \leq \boldsymbol b \label{eq:recons_gen}
\\
&& \boldsymbol x_t \leq \boldsymbol u \notag,
\end{IEEEeqnarray}
where $\rho_t$ is a given time-varying input (e.g., the electricity price) and $\boldsymbol e$ is an all-ones vector of an appropriate size.
In the technical literature, problem~\eqref{eq:recons_gen} is typically referred to as the \emph{reconstruction problem} or the \emph{forward problem} \cite{chan2014generalized, aswani2015inverse}.

Now assume that the matrix of coefficients $\boldsymbol A$ and the right-hand side vector $\boldsymbol u$ are known and that we are able to observe the multivariate time series $\boldsymbol X' = [\boldsymbol x'_1, \ldots, \boldsymbol x'_T]$, which is presumed to be the solution to the reconstruction problem~\eqref{eq:recons_gen} at every time $t$. That is, $\boldsymbol x'_t$ represents the choices actually made by the decision-maker at time $t$. The basic goal of our inverse optimization approach is to infer the unknown parameter vectors $\boldsymbol c$  and $\boldsymbol b$ from $\boldsymbol X'$ given $\boldsymbol A$, $\boldsymbol u$, and the series of measured inputs $\rho_t$. To this end, one tries to find values for the unknowns $\boldsymbol c$  and $\boldsymbol b$ such that the observed choices $\boldsymbol X'$ are as optimal as possible for every problem~\eqref{eq:recons_gen}. With this aim in mind, we solve the following generalized inverse optimization problem:
\begin{IEEEeqnarray}{lrlc}
\!\!\!\!\!\!\!\!\text{GIOP: }& \underset{\boldsymbol b,\boldsymbol \phi_t,\boldsymbol \lambda_t,\boldsymbol c}{\text{Minimize}} \ &   \sum_{t=1}^T \epsilon_t \notag
\\
&\text{subject to} \  &  \boldsymbol b^T \boldsymbol \lambda_t + \boldsymbol u^T \boldsymbol \phi_t  -  \epsilon_t =  (\boldsymbol c-\rho_t \boldsymbol e)^T   \boldsymbol x'_t   \quad  & \!\!\forall t \notag
\\
&& [  \boldsymbol A^T \  \boldsymbol I ]  [  \boldsymbol \lambda^T_t \  \boldsymbol  \phi^T_t ]^T =  (\boldsymbol  c -\rho_t \boldsymbol e) & \!\!\forall t \label{eq:inv}
\\
&&   \boldsymbol A \boldsymbol x'_t \leq \boldsymbol b \notag  & \!\!\forall t
\\
&&  \boldsymbol \phi_t, \boldsymbol \lambda_t, \epsilon_t \geq 0 & \!\!\forall t \notag
\end{IEEEeqnarray}
where $\boldsymbol I$ is the identity matrix of an appropriate size. The objective of optimization problem~\eqref{eq:inv} is to minimize the sum over time of the duality gaps associated with the primal-dual reformulation of problem \eqref{eq:recons_gen}. Thus, when the objective function of GIOP is equal to zero, namely, the accumulated duality gap is zero, $ \boldsymbol x'_t$ is optimal in RP$_t$($\rho_t|$$\boldsymbol c, \boldsymbol b$),  $\forall t$.
The first and second constraints in \eqref{eq:inv} are the relaxed strong duality condition and the dual problem constraints of \eqref{eq:recons_gen}, respectively. The third inequality represents the primal feasibility constraint involving the unknown right-hand side vector $\boldsymbol b$.
The second primal constraint in \eqref{eq:recons_gen} is, in contrast, omitted in~\eqref{eq:inv}, because it does not involve any decision variable in GIOP.

Several challenges arise when solving problem~\eqref{eq:inv}. The most noticeable one is its nonlinear, nonconvex nature, which is the result of the product of variables $\boldsymbol b^T \boldsymbol \lambda_t$ appearing in the strong duality condition. This nonlinearity makes GIOP computationally expensive and hard to solve in general. From this point of view, a method capable of obtaining a good solution to~\eqref{eq:inv} in a reasonable amount of time is needed, even if such a solution may be suboptimal.


Once the parameter vectors $\boldsymbol c$ and $\boldsymbol b$ have been estimated by solving~\eqref{eq:inv}, we can use the reconstruction problem \eqref{eq:recons_gen} to forecast future choices of the decision maker. The estimation of $\boldsymbol c$ and $\boldsymbol b$ through GIOP is anchored in the following two assumptions about the observed choices $ \boldsymbol x'_t$:
\begin{enumerate}
\item Feasibility: $ \boldsymbol x'_t$ is feasible in RP$_t$($\rho_t|$$\boldsymbol c, \boldsymbol b$).
\item Optimality: Given the true $\boldsymbol c$ and $\boldsymbol b$, $ \boldsymbol x'_t$ is optimal in RP$_t(\rho_t)$ and hence, $\epsilon_t$ can be decreased to zero.
\end{enumerate}

These two assumptions, however, do not usually hold in practice for a number of reasons \cite{esfahani2015data}. First, the forward problem \eqref{eq:recons_gen} might be \emph{misspecified} in the sense that it might not represent the actual optimization problem solved by the decision-maker. Therefore, there might not exist $\boldsymbol c$ and $\boldsymbol b$ such that the observed choices $ \boldsymbol x'_t$ are both feasible and optimal for \eqref{eq:recons_gen}. Second, the decision-maker might suffer from \emph{bounded rationality} or \emph{implementation errors}. That is, even if the forward problem \eqref{eq:recons_gen} does prompt the optimal choices to be made by the decision-maker, she might be content with suboptimal choices (due to cognitive or computational limitations, for instance) or there might not be a way to implement such optimal choices without some level of error. Finally, the observed choices $ \boldsymbol X'$ might be corrupted by \emph{measurement noise}.

In this work, though, our intention is to use inverse optimization to \emph{forecast} the future choices of the decision-maker by using the reconstruction problem~\eqref{eq:recons_gen}. This has two important practical implications at least. First, we are not that concerned with the fact that the forward problem \eqref{eq:recons_gen} might be misspecified (this will be indeed the case in the application problem we present later). What we demand from this problem, instead, is that it features good \emph{predictive power} on the futures choices of the decision-maker. In other words, our aim is not to determine values for $\boldsymbol c$ and $\boldsymbol b$ that make the observed choices $ \boldsymbol X'$ both feasible and optimal for \eqref{eq:recons_gen}, but to find the values of these parameters that minimize the \emph{out-of-sample prediction error}.

Given all these practical considerations, in order to compute appropriate values for $\boldsymbol c$ and $\boldsymbol b$, we develop a two-step estimation procedure that deals with the assumptions of feasibility and optimality of $ \boldsymbol X'$ in a \emph{statistical} sense, i.e., with a view to minimizing the out-of-sample prediction error. Furthermore, the proposed two-step estimation procedure overcomes the nonconvexity and computational issues mentioned above regarding the solution to problem~\eqref{eq:inv}.


The first step of the estimation procedure consists in finding a ``good'' feasible region. Note that if $ \boldsymbol b \rightarrow \infty$, then the second constraint in \eqref{eq:recons_gen} is always satisfied. For this reason, we do not want just to find a region for which $\boldsymbol X'$ is  feasible, but the most adequate one in terms of prediction performance. For this purpose, we solve optimization problem~\eqref{eq:feas}, which we refer to as the \emph{feasibility problem} FP(K) and which minimizes a trade-off between the ``infeasibility slack variables'' $\boldsymbol \xi^-$, and the ``feasibility slack variables'' $\boldsymbol \xi^+$, being $0 \leq \text{K} < 1$ a parameter that controls the trade-off between these two quantities.
\begin{IEEEeqnarray}{lrll}
\text{FP($ K$):} \qquad & \underset{\boldsymbol b, \boldsymbol \xi^+, \boldsymbol \xi^-}{\text{Minimize}} \ & \sum_{t=1}^T \Big( \rm{K}  \|\boldsymbol \xi^+_t\| + (1-\rm{K}) \boldsymbol \|\boldsymbol \xi^-_t\| \Big) \notag
\\
&\text{subject to} \  &  \boldsymbol b - \boldsymbol A \boldsymbol  x'_t = \boldsymbol \xi^+_t - \boldsymbol \xi^-_t  &\!\! \forall t \ \  \notag
\\
&&\boldsymbol D^1 \boldsymbol b \leq \boldsymbol d^1  & \label{eq:feas}
\\
&& \boldsymbol \xi^+_t, \boldsymbol \xi^-_t \geq 0 &\!\! \forall t \notag.
\end{IEEEeqnarray}

The value of parameter K is computed by means of cross-validation, as explained below in Section \ref{sec:cv}. Parameters $\boldsymbol D^1$ and $\boldsymbol d^1$ define constraints on $\boldsymbol b$, known a priori, which might be imposed by the nature of $\boldsymbol b$. An example of such a-priori constraint may simply be $\boldsymbol b \geq 0$.

In the second step, we consider $\boldsymbol b$ given as the solution to FP(K), denoted by $\widehat{\boldsymbol b}$. Also, we adjust the observed quantity $\boldsymbol A \tilde{\boldsymbol x}'_t = \boldsymbol A \boldsymbol x'_t - \boldsymbol  \xi^{-*}_t$, where $\boldsymbol \xi^{-*}_t$ is taken from the solution of FP(K). This modification makes $\tilde{\boldsymbol x}'_t$ feasible in RP$_t$($\rho_t|$$\boldsymbol c, \widehat{\boldsymbol b}$). As we show later in Section~\ref{sec:util_estim}, we may need to impose further constraints on the adjusted quantity $\tilde{\boldsymbol x}'_t$ in those cases where it is not univocally determined by  $\boldsymbol A \tilde{\boldsymbol x}'_t = \boldsymbol A \boldsymbol x'_t - \boldsymbol  \xi^{-*}_t$. Then, we solve the following linear programming problem:
\begin{IEEEeqnarray}{lrlc}
\!\!\!\!\!\!\!\text{OP($\widehat{\boldsymbol b}$): } & \underset{ \boldsymbol \phi_t, \boldsymbol \lambda_t, \boldsymbol c}{\text{Minimize}} \ &  \sum_{t=1}^T  \epsilon_t \notag
\\
&\text{subject to} \  &    \widehat{\boldsymbol b}^T \boldsymbol  \lambda_t + \boldsymbol u^T \boldsymbol \phi_t  -   \epsilon_t =  (\boldsymbol c - \rho_t \boldsymbol e)^T   \tilde{\boldsymbol x}'_t  &\ \forall t  \notag
\\
&& [  \boldsymbol A^T \  \boldsymbol I ]  [  \boldsymbol \lambda^T_t \  \boldsymbol  \phi^T_t ]^T =  (\boldsymbol  c - \rho_t \boldsymbol e) &\ \forall t \label{eq:OP}
\\
&&\boldsymbol D^2 \boldsymbol c \leq \boldsymbol d^2  \notag
\\
&&  \boldsymbol \phi_t, \boldsymbol \lambda_t,   \epsilon_t \geq 0 &\ \forall t \notag.
\end{IEEEeqnarray}

The first, the second, and the last constraints in \eqref{eq:OP} are analogue to the ones in \eqref{eq:inv}. The third constraint defines a-priori conditions on $\boldsymbol c$, specified by the parameters $\boldsymbol D^2$ and $\boldsymbol d^2$. The outcome of this problem is the estimated value of the coefficient vector $\boldsymbol c$, named as $\widehat{ \boldsymbol c}$.

Finally, given $\widehat{ \boldsymbol b}$, $\widehat{ \boldsymbol c}$ and $\rho_{T+1}$, we can forecast the future decision-maker's choices by solving RP$_{T+1}$($\rho_{T+1}|$$\widehat{ \boldsymbol c},\widehat{ \boldsymbol b}$).

\subsection{Statistical Determination of K} \label{sec:cv}

In practice, we find the value of the penalty parameter K using cross-validation \cite[Ch. 7,]{friedman2001elements}. We partition $\boldsymbol X'$ in three subsets: the training set $\boldsymbol X'_{tr}$, the validation set $\boldsymbol X'_{val}$, and the test set.
In a few words, for each given value of K, we use the training set for parameter fitting and the validation set to asses the forecasting performance. The best choice of K is, thus, the one that minimizes the out-of-sample prediction error.


The advantage of using this approach is threefold. First, by tuning the value of K for the validation set, we seek to minimize the out-of-sample prediction error (a criterion specially suited for forecasting purposes). Second, we solve three LP problems for each tested value of K, hence, the computational burden of the tuning algorithm is relatively low. Finally, the evaluations of different values of K are independent of each other so they can be executed in parallel.

\subsection{Leveraging Auxiliary Information}
When using inverse optimization for \emph{forecasting} a time series, it is relevant to consider the case where we let the unknown parameter vectors $\boldsymbol c$ and $\boldsymbol b$ vary over time so as to capture structural changes in the decision-making problem~\eqref{eq:recons_gen}. To this end, we assume that we also observe  a number of time-varying \emph{regressors} $\boldsymbol Z_t$ that, to a lesser or greater extent, may affect the decision maker's choices. We then describe the unknown vectors $\boldsymbol c$ and $\boldsymbol b$ as functions of those regressors by letting $\boldsymbol c_t = f_c(\boldsymbol Z_t)$ and $\boldsymbol b_t = f_b(\boldsymbol Z_t)$ in problems \eqref{eq:feas} and \eqref{eq:OP}, respectively. In this way, functions $f_c(\cdot)$ and $f_b(\cdot)$ become decision variables in our estimation problem. In this paper, we consider $f_c(\cdot)$ and $f_b(\cdot)$ to be affine functions. Hence, the inverse optimization problem seeks the most optimal set of intercepts and affine coefficients that relate $\boldsymbol Z_t$ with $\boldsymbol c_t$ and $\boldsymbol b_t$, as we exemplify below.
Note that the past choices of the decision-maker, namely, $\boldsymbol X'$, can also be treated as regressors.


\section{Methodology Applied to Forecast Price-responsive Loads} \label{sec:applied}
%

In this section we illustrate the use of the proposed inverse optimization approach to forecast the aggregate power load of a pool of price-responsive consumers.


We consider that the available information is the measured power consumption $x'_t$ of the pool; the electricity price $p_t$, which is broadcast to every load in the pool; and the realizations of a set of explanatory variables $Z_{r,t}$ for every time period $t$. The aggregate response $\boldsymbol x_t$ of the loads to the price of electricity at time $t$ is assumed to be the solution to the following forward/reconstruction problem:
\begin{subequations} \label{eq:RPa}
\begin{IEEEeqnarray}{llr}
\underset{ \boldsymbol x_t}{\text{Maximize}} \ & \sum_{b=1}^B x_{b,t} \left(u_{b,t} - p_{t} \right) \label{eq:RPa_obj}
\\
\text{subject to} &\underline{P}_t \leq \sum_{b=1}^B x_{b,t} \leq \overline{P}_t &(\underline{\lambda}_{t},\overline{\lambda}_{t})
\label{eq:RPa_Pbound}
\\
 & 0 \leq  x_{b,t} \leq E_{b,t}  & (\underline{\phi}_{b,t},\overline{\phi}_{b,t}) \ \ \forall b \label{eq:RPa_Ebound}.
\end{IEEEeqnarray}
\end{subequations}

Problem~\eqref{eq:RPa} takes the form of \eqref{eq:recons_gen}. Its objective function \eqref{eq:RPa_obj}, to be maximized, represents the aggregate consumers' surplus or welfare, given as the product of the pool consumption and the difference between the marginal utility and the electricity price. We consider a step-wise marginal utility curve made up of $B$ blocks, each of a width $E_{b,t}$, as enforced by  \eqref{eq:RPa_Ebound}, and a value $u_{b,t}$. The aggregate load of the pool, given as $\sum_{b=1}^B x_{b,t}$, is bounded from below and above by $\underline{P}_t$ and $\overline{P}_t$, respectively, as expressed by \eqref{eq:RPa_Pbound}. Symbols within parentheses correspond to the dual variables associated with each constraint.


The goal of our inverse optimization methodology is to estimate appropriate values for $u_{b,t}$, $\underline{P}_t$ and $\overline{P}_t$, based on the observed $x'_t$, $p_t$, and $Z_{r,t}$, such that  the solution to the reconstruction problem \eqref{eq:RPa} serves as a good forecast of the future aggregate power consumption $x_{t+1}$ of the pool of loads. For this purpose, we employ the estimation procedure outlined in Section \ref{sec:method}.
Note that the width of each block, $E_{b,t}$, is treated as a parameter and need not to be estimated. Later, in Section \ref{sec:util_estim}, we give a practical rule for fixing it.

Next we provide concrete formulations for the estimation of $\underline{P}_t, \overline{P}_t$ and $u_{b,t}$. The problem of estimating the bounds $\underline{P}_t$ and $\overline{P}_t$, which we refer to as the \emph{bound estimation problem}, is presented in Section \ref{sec:bound_estim}. The problem of estimating the marginal utilities $u_{b,t}$, which we call  the \emph{marginal utility estimation problem}, is presented in Section \ref{sec:util_estim}. A discussion about the proposed methodology is finally given in \ref{sec:discussion}.

\subsection{Bound Estimation Problem} \label{sec:bound_estim}
The bound estimation problem is derived from the feasibility problem \eqref{eq:feas} and consists in determining the bounds $\underline{P}_t$ and $\overline{P}_t$ by minimizing the following objective function:
\begin{subequations} \label{eq:FPa}
\begin{IEEEeqnarray}{ll}
\underset{ \underline{\boldsymbol P}, \overline{\boldsymbol P}, \boldsymbol \xi,\boldsymbol \mu, \boldsymbol \alpha }{\text{Minimize}} \  \sum_{t=1}^T  \Bigg( & (1-\rm{K}) \left(  \overline{\xi}^{+}_t + \underline{\xi}^{+}_t \right) +
\notag \\
 &\rm{K} \left( \overline{\xi}^{-}_t + \underline{\xi}^{-}_t\right) \Bigg) \label{eq:FPa_Obj}
\end{IEEEeqnarray}
subject to
\begin{IEEEeqnarray}{lr}
\overline{P}_t - x_t' =  \overline{\xi}^{+}_t - \overline{\xi}^{-}_t   \qquad \qquad & \forall t
\\
x_{t}' - \underline{P}_t =  \underline{\xi}^{+}_t - \underline{\xi}^{-}_t   & \forall t
\\
\underline{P}_{t} \leq \overline{P}_{t} & \forall t \label{eq:FPa_MaxMin}
\\
\underline{P}_t = \underline{\mu} + \sum_{r=1}^R \underline{\alpha}_i Z_{r,t}  & \forall t  \label{eq:FPa_Pmin}
\\
\overline{P}_t = \overline{\mu} + \sum_{r=1}^R \overline{\alpha}_i Z_{r,t}    & \forall t \label{eq:FPa_Pmax}
\\
0 \leq \overline{\xi}^+_t, \overline{\xi}^-_t,\underline{\xi}^+_t, \underline{\xi}^-_t  & \forall t
\end{IEEEeqnarray}
\end{subequations}
where $\boldsymbol \xi = \left[\overline{\boldsymbol \xi}^{+}; \underline{\boldsymbol \xi}^{+}; \overline{\boldsymbol \xi}^{-}; \underline{\boldsymbol \xi}^{-}\right]$, $\boldsymbol \mu = \left[\overline{\boldsymbol \mu}; \underline{\boldsymbol \mu}\right]$, and $\boldsymbol \alpha = \left[\overline{\boldsymbol \alpha}; \underline{\boldsymbol \alpha}\right]$.

The objective function \eqref{eq:FPa_Obj} comprises two terms, weighted by the parameter K with $0 \leq \rm{K}<1$. The first and second terms represent the amount of measured load that falls inside and outside the interval $[\underline{P}_t, \overline{P}_t]$, respectively.

Constraint \eqref{eq:FPa_MaxMin} ensures that the estimated lower bound is always lower than the upper bound.
Constraints \eqref{eq:FPa_Pmin} and \eqref{eq:FPa_Pmax} impose an affine relationship between the regressors and the load bounds. We denote the estimates of the lower and upper bounds at the optimum as $\widehat{\underline{P}}_t$ and $\widehat{\overline{P}}_t$, respectively.
It is worth mentioning the special case where previous load observations $x'_{t-1} \ldots x'_{t-l}$ are included as regressors, in a similar way as traditional auto-regressive models do \cite{MadsenBook}. Also, it should be noted that we do not treat the price at time $t$ as a regressor here, since its effect is captured through the objective function of the forward problem by solving the optimality problem \eqref{eq:OPa}.

The parameter K is computed using the cross-validation approach outlined in Section \ref{sec:cv}. Values of K close to 1 yield a ``wide'' interval $[\widehat{\underline{P}}_t, \widehat{\overline{P}}_t]$, whereas values of K close to zero produce a narrow interval. Therefore, K can be interpreted as an indicator of the \emph{price responsiveness} of the load, since the precise value that the load will take on within the interval $[\widehat{\underline{P}}_t, \widehat{\overline{P}}_t]$ is left to be explained by the electricity price. Notice that when $\text{K}= 0$ the bound estimation problem boils down to fitting an ARX model by minimizing the MAE, because in this case it holds that $\widehat{\underline{P}}_t= \widehat{\overline{P}}_t$.

Now consider the solution to \eqref{eq:FPa}. In order for the reconstruction problem \eqref{eq:RPa} to be feasible, the estimated bounds must satisfy that $\widehat{\underline{P}}_{t} \leq \widehat{\overline{P}}_{t}$. This is enforced for the training data set by Equation \eqref{eq:FPa_MaxMin}, but it is not necessarily satisfied for any \emph{plausible} data point outside this set. Generally, when forecasting, consistent bounds $\underline{P}_{t} \leq \overline{P}_{t}$ must be obtained for all possible future realizations of the regressors $Z_{r,t}$. We can ensure this by robustification \cite{ben2009robust}, as done in \cite{Saez2015}.

\subsection{Marginal Utility Estimation Problem} \label{sec:util_estim}

The marginal utility estimation problem is derived from the optimality problem \eqref{eq:OP} once the bounds $\widehat{\underline{P}}_t$ and $\widehat{\overline{P}}_t$ have been estimated by solving problem \eqref{eq:FPa}.

Prior to solving the marginal utility estimation problem, the measured load is adjusted so that it becomes feasible in the reconstruction problem \eqref{eq:RPa}. For this purpose, we define the adjusted load as $\tilde{x}'_t =  x'_t - \overline{\xi}^{-*}_t + \underline{\xi}^{-*}_t$, where $\overline{\xi}^{-*}_t$  and $\underline{\xi}^{-*}_t$ are taken from the solution to problem \eqref{eq:FPa}. Note that this is equivalent to defining $\tilde{x}'_t$ as
\begin{IEEEeqnarray}{cl}
\tilde{x}'_t = \sum_{b=1}^B \tilde{x}'_{b,t} =
 \begin{cases}
     \widehat{\underline{P}}_t       & \quad \text{if } x'_t < \widehat{\underline{P}}_t
     \\
    x'_t  & \quad \text{if }  \widehat{\underline{P}}_{t} \leq x'_t \leq \widehat{\overline{P}}_{t} \\
   \widehat{\overline{P}}_{t}  & \quad \text{if } x'_t > \widehat{\overline{P}}_{t} \\
 \end{cases} \label{eq:split}
\end{IEEEeqnarray}
with $\tilde{x}'_{b,t} \leq E_{b,t}$. We further impose that the load blocks are to be filled in sequential order starting with $b =1$.

To fix the width of each load block, we proceed as follows. We set the width of the first block to be equal to the lower bound, namely, $E_{1,t} =\widehat{\underline{P}}_t$.
The width of the remaining blocks is computed such that the interval $[\widehat{\underline{P}}_{t}, \widehat{\overline{P}}_{t}]$ is equally divided, that is, $E_{b,t} = (\widehat{\overline{P}_t} - \widehat{\underline{P}_t})/(B-1), \ \forall b>1, \forall t$. In order for the lower bound to be effective, we set the marginal utility for the first block (denoted as $u_{1,t}$) to a large number. This is done in Equation \eqref{eq:util_maxu}. Consequently, at the optimum, the first block of energy is always filled with $x_{1,t} = E_{1,t}=\widehat{\underline{P}}_t$, since its corresponding marginal utility is always higher than the electricity price. This practical rule allows us to enforce the lower bound through the use of $E_{1,t}$ and $u_{1,t}$.

The optimization variables in the marginal utility estimation problem \eqref{eq:OPa} are $\boldsymbol \Omega = \lbrace \boldsymbol \epsilon,  \boldsymbol u , \boldsymbol \mu^u, \boldsymbol \alpha^u,  \underline{\boldsymbol \lambda}, \overline{\boldsymbol \lambda},  \underline{\boldsymbol \phi},  \overline{\boldsymbol \phi} \rbrace$. This problem aims to minimize the sum of duality gaps $\epsilon_t$ of the reconstruction problem \eqref{eq:RPa}, that is, to find the marginal utilities $u_{b,t}$ such that the adjusted observed load $\tilde{x}'$ is as optimal as possible.

{\allowdisplaybreaks
\begin{subequations} \label{eq:OPa}
\begin{IEEEeqnarray}{rlc}
 \underset{ \boldsymbol \Omega }{\text{Minimize}} &\  \sum_{t=1}^T \epsilon_t &
\\
\text{subject to} &  \  \widehat{\overline{P}}_t  \overline{\lambda}_{t} - \widehat{\underline{P}}_t  \underline{\lambda}_{t} + \sum_{b=1}^B  E_b \overline{\phi}_{b,t} - \epsilon_t
= &\notag
\\
& \qquad \sum_{b=1}^B \tilde{x}_{b,t}' \left(u_{b,t} - p_{t} \right)  \qquad  & \forall t \label{eq:util_strong}
\\
& - \underline{\phi}_{b,t} + \overline{\phi}_{b,t} -\underline{\lambda}_{t} + \overline{\lambda}_{t}  = u_{b,t} - p_{t}     & \forall b ,t \label{eq:util_stat}
\\
& u_{b,t} = \mu^u_b + \sum_r \alpha_r^u Z_{r,t}	 & \forall b , t \label{eq:util_defu}
\\
& \mu^u_b \geq \mu^u_{b+1}   & \forall b<B \label{eq:util_decre}
\\
& \mu^u_1 \geq 200+\mu^u_{2} \label{eq:util_maxu}
\\
& 0 \leq \overline{\lambda}_{t},\underline{\lambda}_{t}, \underline{\phi}_{b,t},\overline{\phi}_{b,t}	& \forall b, t \label{eq:util_pos1}.
\end{IEEEeqnarray}
\end{subequations}
}

Constraint \eqref{eq:util_strong} defines the relaxed strong duality condition, with the objective function of the dual of problem~\eqref{eq:RPa} minus the duality gap on the left-hand side, and the primal objective function on the right-hand one.
Equations \eqref{eq:util_stat} are the constraints of the dual of problem \eqref{eq:RPa}. Constraint \eqref{eq:util_defu} defines the marginal utilities as affine combinations of the regressors. Constraint \eqref{eq:util_decre} forces the estimated marginal utility to be monotonically decreasing, and constraint \eqref{eq:util_maxu} imposes a high utility for the first block. Finally, constraint \eqref{eq:util_pos1} enforces the non-negative character of the dual variables.


\subsection{Discussion} \label{sec:discussion}

The proposed inverse optimization framework can be seen as a generalization of a linear time series model: the relationship between the load and the regressors is linear, but the relationship between the load and the price at time $t$ is not.

Recall that the proposed methodology is composed by two problems that are solved sequentially: the feasibility problem \eqref{eq:FPa} and the optimality problem \eqref{eq:OPa}. In the feasibility problem, we model the linear relationship between the load and the regressors, excluding the price at time $t$. The penalty parameter K is optimally chosen by cross-validation, and it affects the width of the interval $[\widehat{\underline{P}}_t, \widehat{\overline{P}}_t]$.
Afterwards, in the optimality problem we model the non-linear relationship between the load that falls inside the interval $[\widehat{\underline{P}}_t, \widehat{\overline{P}}_t]$, and the price at time $t$.
For this reason, a narrow interval implies that the variability of the load left to be explained by the price at time $t$ is very small. On the other hand, a wide interval indicates that the load can be explained by the price at time $t$ to a large extent. Its non-linear relationship is estimated by the optimality problem. Unlike in the proposed scheme, in a simple linear regression model, the relationship between the load and the price is given by an affine coefficient.

\section{Simulation of Price-responsive Buildings} \label{sec:simulation}

We simulate the price-response behavior of a pool of buildings equipped with heat pumps. To this end, we  use the work in \cite{Halvgaard2012}. The heating dynamics of each building is described by a state-space model that consists of three states: indoor air temperature $y^r_{t}$, floor temperature $y^f_{t}$, and temperature of the water $y^w_{t}$ inside a tank connected to a heat pump. The only input is the electricity consumption $x_t$. The state-space model writes as follows, where $\b y_t = [ y^r_{t}, y^f_{t}, y^w_{t} ]^{T}$:
\begin{IEEEeqnarray}{llr}
\b y_{t} = \text{\bf A} \b y_{t-1} +  \text{\bf B}  x_{t-1} + \text{\bf E}  \b z_{t-1} 	& \qquad \forall t 	\label{eq:state}
\end{IEEEeqnarray}
where \textbf{A}, \textbf{B} and \textbf{E} are the matrices of the coefficients defining the state-space model. The temperature of the air outside the building $z_t^a$ and the solar irradiance $z_t^s$ are considered as external disturbances in $\b z_t  = [z_t^a, z_t^s]^{T}$.

The heat pump in each building schedules its consumption by solving an Economic Model Predictive Control (EMPC) problem that minimizes the cost of its consumption plus a penalty term for not complying with a comfort band:
\begin{subequations}\label{eq:EMPC}
\begin{IEEEeqnarray}{llc}
  \underset{ \b y, \b x ,\b v }{\text{Minimize}} \ & \sum_{t=1}^T p_t  x_t + \rho v_t  	\label{eq:sim-obj}
\\
\text{subject to} \ \ \ & \b y_{t} = \text{\bf A} \b y_{t-1} + \text{\bf B} x_{t-1} + \text{\bf E}  \b z_{t-1} 	\ &  \forall t 	\label{eq:sim-state}
\\
& 0 \leq x_t \leq x^{max} 	& \forall t 	\label{eq:sim-minmax}
\\
& y_{t,min}^r \leq y_t^r + v_t   	& \forall t 	\label{eq:sim-pen1}
\\
& y_{t,max}^r \geq y_t^r - v_t   	& \forall t 	\label{eq:sim-pen2}
\\
& v_t \geq 0    	& \forall t 	\label{eq:sim-posit}.
\end{IEEEeqnarray}
\end{subequations}

The objective function \eqref{eq:sim-obj} minimizes the cost of purchasing $x_t$ kWh of energy at the price $p_t$, with a penalization of $\rho v_t$ if the room temperature is not within the desired comfort band. Equation \eqref{eq:sim-state} determines the time evolution of the states of the model. The maximum power consumption of the heat pump, set by Equation \eqref{eq:sim-minmax}, is $x^{max}$ kW. Finally, equations \eqref{eq:sim-pen1}, \eqref{eq:sim-pen2} and \eqref{eq:sim-posit} define the comfort temperature band, given by $y_{t,min}^r$ and $y_{t,max}^r$, and the slack variable $v_t$.

\begin{figure}
  \centering
  \includegraphics[width=.5\textwidth]{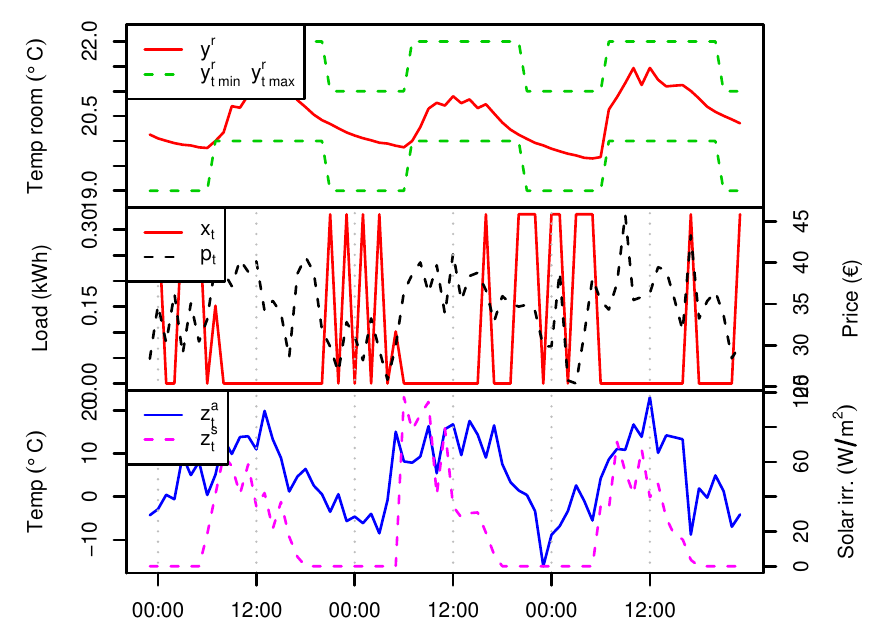}
  \caption{Evolution of room temperature, comfort bands, price, load, and disturbances for one building during 72 hours}
  \label{fig:HP}
\end{figure}

The values for the coefficients of the EMPC model \eqref{eq:EMPC} are taken from \cite[Tab. 1]{Zugno2013}. The effect of the solar irradiance on the room temperature is set to be equal to $\SI{0.01}{\degreeCelsius} /(W/m^2)$, and of $\SI{0.001}{\degreeCelsius} /(W/m^2)$ on the floor temperature.
An example of the behavior of one building during 72 hours, with hourly observations, is depicted in Fig.~\ref{fig:HP}.

A total of 100 buildings are simulated by randomly perturbing the heat-transfer coefficients that define matrix \text{\bf A}.
Modifying slightly these coefficients allows us to simulate the behavior of buildings with different structural characteristics. The perturbations are randomly drawn from a uniform distribution centered around zero with a variance equal to 1/50 the magnitude of the corresponding coefficient. The magnitude of the perturbations is chosen high enough so that different building structures are modeled, but not too high so that the state-space system becomes unstable. The magnitude of the perturbations has been chosen by trial-and-error, and its effectiveness is proven to be useful as explained in the remaining of this section and in the case study of Section \ref{sec:study_case}.

\begin{figure}
  \centering
  \includegraphics{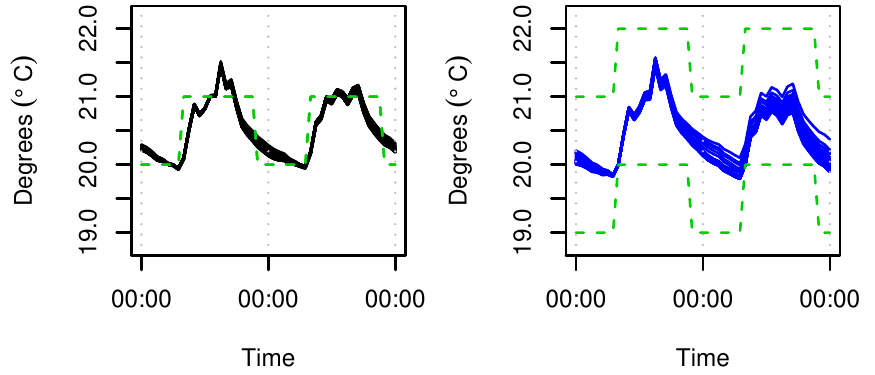}
  \caption{Room temperature and comfort bands. Left: no-flex case; right: flex case}
  \label{fig:Tr}
\end{figure}

We simulate the behavior of two classes of buildings and aggregate the simulated information in two data sets. In the first one, called \textit{no flex}, the comfort bands for the temperature inside the room are equal to each other. In the second case, called \textit{flex}, the comfort bands for the temperature of the air inside the room are \SI{2}{\degreeCelsius} apart from each other. A sample of the simulated comfort bands and temperatures inside the rooms is shown in Fig. \ref{fig:Tr}. On the left plot, we show the \textit{no flex} case. Naturally, the temperature inside the room is as close as possible to the desired one. On the right plot, for the \textit{flex} case, the temperature inside the room features a higher variation across buildings.

\begin{figure}
  \centering
  \includegraphics{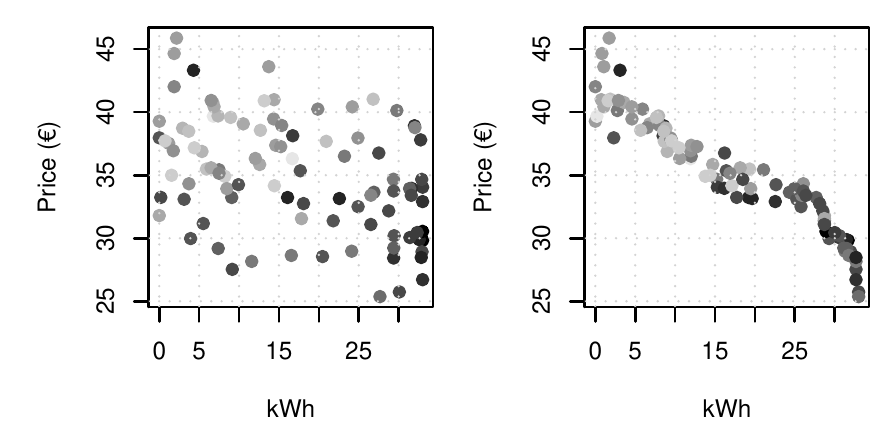}
  \caption{Price and aggregate load for the non-flexible (left) and flexible (right) cluster of buildings. Dark tonalities indicate low ambient temperature.}
  \label{fig:PriceVSload}
\end{figure}

The effect of the electricity price on the aggregated load, for the two data sets, is displayed in Fig. \ref{fig:PriceVSload}. On the left plot, the \textit{no flex} data set shows barely no relationship between load and price. On the other hand, on the right plot, the \textit{flex} data set shows a clear non-linear relationship. In both plots, black colors indicate that the temperature of the outside air is low. Naturally, the aggregated load is higher at times with low ambient temperature because of the need for heating up the water tank.

To sum up, the simulated data sets seem to be fair representations of the behavior of a pool of price-responsive buildings, hence, we proceed to use the simulated data sets for testing the performance of the proposed load forecasting method.

\section{Case Study} \label{sec:study_case}

We now asses the performance of the reconstruction problem \eqref{eq:RPa} when forecasting the load of the pool of buildings one step ahead, that is, one hour in advance. In problem \eqref{eq:RPa}, the marginal utilities, the minimum power, and the maximum power are calculated using the methodology introduced in Section \ref{sec:method} and \ref{sec:applied}.  We compare the performance of the proposed methodology using the \textit{flex} and \textit{no flex} data sets, which were introduced in Section \ref{sec:simulation}.

%
%

The regressors that we consider for describing the dynamics of the estimated parameters in \eqref{eq:FPa_Pmin}, \eqref{eq:FPa_Pmax} and \eqref{eq:util_defu}, are the hour indicator, outside temperature, solar irradiance, and historical lagged price and load data.
At every time period $t$, we assume that the regressors up to that period are known. We also assume that the price at time $t$ is known. The training set consists of 505 data points, that is, three weeks of data. Furthermore, we set the total number of load blocks to $B=20$.

\begin{figure}
  \centering
  \includegraphics{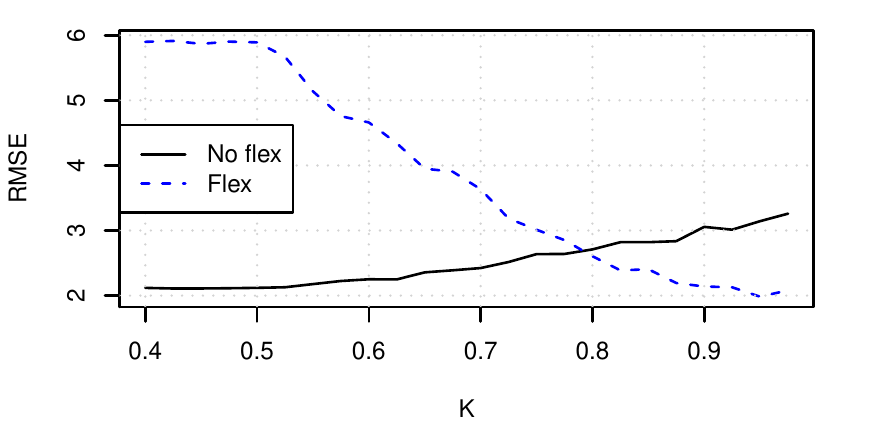}
  \caption{The RMSE of the 1-steap ahead predictions is shown for different values of K, using the \textit{no flex} dataset (dashed) and the \textit{flex} dataset (continuous)}
  \label{fig:CV}
\end{figure}

The first step in the estimation procedure is to tune parameter K in \eqref{eq:FPa} following the cross-validation strategy from Section \ref{sec:cv}. The results, displayed in Fig. \ref{fig:CV}, show the Root Mean Square test Error (RMSE) for the two data sets, using a validation period of one week of data. For the \textit{flex} dataset, the optimal value of K turns out to be 0.98. Recall that values of K close to 1 indicate that the interval $[\widehat{\underline{P}},\widehat{\overline{P}}]$ is wide. Therefore, in the considered application, this also means higher responsiveness of the load to the price.
The continuous line in Fig. \ref{fig:CV} represents the RMSE for the \textit{no flex} dataset, for different values of K. The best forecasting performance is achieved for $\text{K} \leq 0.5$. In the \textit{no flex} case, the load is independent of the price. Consequently, the best forecasts are achieved when the interval $[\widehat{\underline{P}},\widehat{\overline{P}}]$ is very small, namely, when $\widehat{\underline{P}}=\widehat{\overline{P}}$. It is noteworthy to say that the solution in this case is equivalent to fitting an autoregressive linear model with exogenous inputs by minimizing the mean absolute error.

On average, when using 3 weeks of data, 20 load blocks, and 38 regressors, the time for the whole estimation process takes around 10 seconds on a personal Linux-based machine with 4 cores clocking at 2.90GHz and 12 GB of RAM. R and CPLEX 12.3 under GAMS are used to process the data and solve the optimization models. We conclude that, because of its low computational requirements, the proposed methodology is attractive for implementation in a real-life setup.

\subsection{Benchmark on a Test Period} \label{sec:bench}

We benchmark the forecasting capability of the proposed methodology against two other methods. The first one is a simple persistence model, where the forecast load at time $t$ is set to be equal to the observed load at $t-1$.
The second model is an Autoregressive Moving Average Model with eXogneous inputs (ARMAX)  \cite[Ch. 5]{MadsenBook}, similar to the one used in \cite{corradi2012}. The aggregate load $x_{t}$ is modeled as a linear combination of the past values of load, past errors, and regressors. In mathematical terms, the ARMAX model can be written as
\begin{IEEEeqnarray}{l}
x_t = \mu + \epsilon_t +  \sum_{p=1}^P \varphi_p x_{t-p} + \sum_{r=1}^R \gamma_r Z_{r} + \sum_{q=1}^Q \theta_q \varepsilon_{t-q}
\end{IEEEeqnarray}

\noindent with $\epsilon_{t} \sim \text{N(0, }\sigma^2 \text{)}$ and $\sigma^2$ being the variance. The optimal combination of $P$ and $Q$ is chosen according to the AICc criteria \cite{burnham2004multimodel}. In order to make reasonable comparisons, the same explanatory variables are used for our inverse-optimization-based model and for the ARMAX, including the price at time $t$. Recall that the price at time $t$ is considered in the optimality problem \eqref{eq:OPa} but not in the feasibility problem \eqref{eq:FPa}.

\begin{figure}
  \centering
  \includegraphics{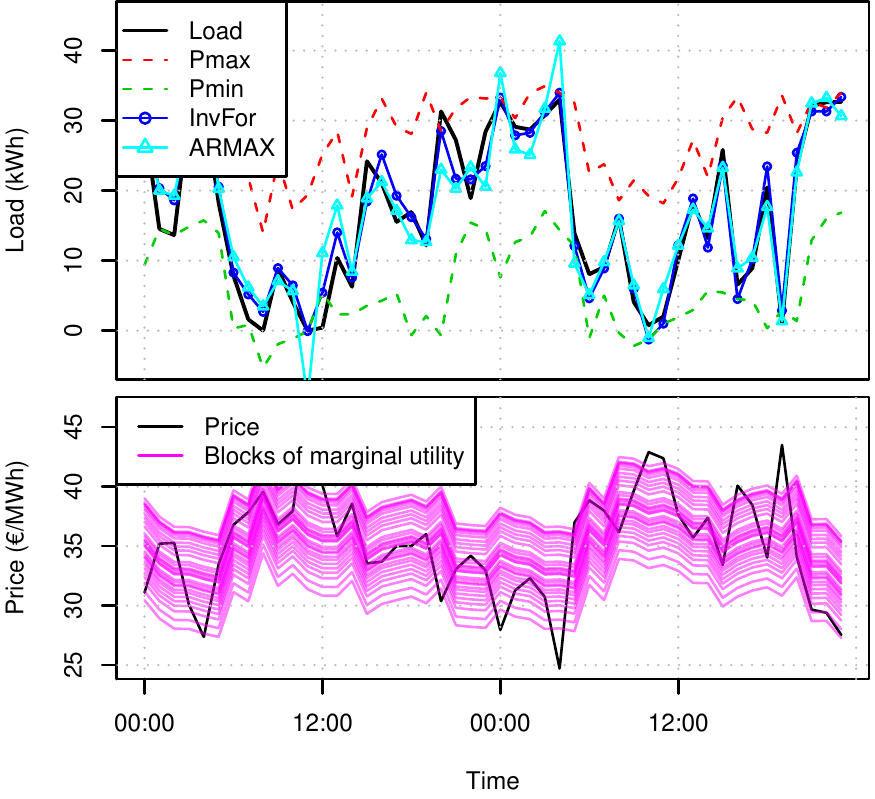}
  \caption{On the top, the actual load is displayed together with the predictions from the inverse-optimization methodology and the ARMAX model. On the bottom, the price is shown together with the estimated marginal utility blocks}
  \label{fig:Test_pred}
\end{figure}

We run 1-step ahead predictions in a rolling-horizon manner for a period of 5 days, re-estimating the parameters at every hour. The upper plot of Fig.~\ref{fig:Test_pred} shows the actual aggregate load together with those predicted by  the proposed inverse-optimization model ({\it InvFor}) and the ARMAX. The estimated minimum and maximum load bounds are able to explain a certain part of the variability of the load. The remaining variability is explained by the relationship between the marginal utilities and the price. The predictions made by the ARMAX are also able to anticipate the behavior of the load, but to a lesser extent. On the bottom plot of Fig. \ref{fig:Test_pred}, the electricity price is displayed together with the estimated marginal utility blocks, for each hour of the test period. The magnitude and distribution of the marginal utilities change with time and capture the dynamic response of the load to the price.

Performance metrics computed over the test set are summarized in Table \ref{tab:benchmark}. Each row is relative to one of the three benchmark models. Columns 2 and 4 give information on the Normalized Root Mean Square Error (NRMSE), defined as
\begin{IEEEeqnarray}{l}
{\rm NRMSE} =\frac{1}{x^{max} - x^{min}}\sqrt{ \frac{1}{T} \sum_{t=1}^T\left(\sum_{b=1}^B\widehat{x}_{b,t} - x'_t\right)}
\end{IEEEeqnarray}
\noindent and columns 3 and 5 on the Symmetric Mean Absolute Percentage Error (SMAPE)
\begin{IEEEeqnarray}{l}
{\rm SMAPE} =\frac{1}{T} \sum_{t=1}^T \frac{|\sum_{b=1}^B\widehat{x}_{b,t} - x'_t|}{(|\sum_{b=1}^B\widehat{x}_{b,t}| - |x'_t|)/2}.
\end{IEEEeqnarray}

In Table \ref{tab:benchmark} we also compare the performance of the proposed forecasting method using the two simulated data sets. On the left part, we show the performance measures relative to the \textit{no flex} data set. The ARMAX and the \textit{InvFor} models yield almost identical results in terms of NRMSE and SMAPE. This is indeed reasonable because, as mentioned earlier in Section \ref{sec:method}, the \textit{InvFor} model with a penalty parameter of $K=0$ is equivalent to fitting an ARX.

The differences between the ARMAX and the \textit{InvFor} stand out when used for predicting the \textit{flex} data set. On the right side of Table \ref{tab:benchmark}, we see that our methodology outperforms the ARMAX with a NRMSE and a SMAPE 32\% and 16.8\% lower, respectively. The persistence model, as expected, exhibits the worst performance. We conclude that the non-linear relationship between the price and the load is well captured by the \textit{InvFor} model.

\begin{table}[ht]
\centering
\caption{Benchmark for the test set}
\begin{tabular}{ccccc}
& \multicolumn{2}{c}{ No Flex} & \multicolumn{2}{c}{ Flex}  \\
  & NRMSE & SMAPE  & NRMSE & SMAPE\\
  \midrule
  \textit{Persistence}  & 0.1727 & 0.1509  & 0.3107 & - \\
  \textit{ARMAX}      & 0.10086 & 0.08752  & 0.13107 & 0.08426 \\
  \textit{InvFor}       & 0.10093 & 0.0886  & 0.08903 & 0.07003  \\
   \hline
\end{tabular}
\label{tab:benchmark}
\end{table}


\section{Conclusion} \label{sec:conclusions}

This paper proposes a new method to forecast price-responsive electricity consumption. The price response is described by an optimization problem, which is characterized by a set of unknown parameters. The problem of estimating these parameters is nonlinear and nonconvex. We formulate a two-step algorithm to statistically approximate its solution, where in each step we solve a linear problem.
The proposed approach is data-driven and makes use of a cross-validation scheme to minimize the out-of-sample prediction error. Moreover, a set of regressors is used to explain the variability of the price-response of the load.

A simulation framework is used to asses the performance of the proposed methodology. The simulation comprises a set of price-responsive buildings equipped with a heat pump. The presented methodology is used for 1-step ahead predictions. Results show that the non-linear relationship between the price and the aggregate load is successfully captured and that the proposed method outperforms well-known benchmark models.


\bibliographystyle{IEEEtran}


\end{document}